\documentclass[12pt,letterpaper]{article}
\usepackage[utf8]{inputenc}
\usepackage[T1]{fontenc}
\usepackage{amsmath, amssymb, amsthm}
\usepackage{geometry}
\usepackage{hyperref}
\theoremstyle{plain}
\newtheorem{theorem}{Theorem}[section]
\newtheorem{proposition}[theorem]{Proposition}
\newtheorem{lemma}[theorem]{Lemma}
\newtheorem{corollary}[theorem]{Corollary}
\usepackage{authblk} 
\theoremstyle{definition}
\newtheorem{definition}[theorem]{Definition}
\newtheorem{remark}[theorem]{Remark}

\usepackage{abstract}

\title{A cubic form on a cubic threefold-hyperplane pair}

\author{Gustavo Serrano Herrera}
\date{}

\begin{document}	
	\maketitle
	\begin{abstract}
		In this paper, we construct a non-degenerate cubic form for a specific class of pairs $(X,Y)$, where $X$ is a smooth cubic threefold and $Y$ is a hyperplane section. By establishing explicit isomorphisms between twisted logarithmic cohomology groups and the graded pieces of an Artinian-Gorenstein algebra, the logarithmic Jacobian ring, we study abstract cohomological data with concrete computational algebraic methods, extending the classical Yukawa-Griffiths constructions to this setting.
	\end{abstract}
	
	\section*{Introduction}
	The geometry of compact projective hypersufaces have been deeply studied by the work of Phillip Griffiths. In \cite{Griffiths}, Griffiths established a correspondence between the primitive cohomology groups of a smooth projective hypersurface and the graded components of its associated Jacobian ring. In that classical framework, Macaulay's duality theorem provides a perfect pairing on the Jacobian ring and allows the explicit computation of the Yukawa-Griffiths cubic form via elementary polynomial multiplication. This cubic form defined through the cup product of infinitesimal deformations and period integrations, measures the variations of complex structures, serving as the algebraic engine behind the infinitesimal Torelli theorem and encode other local geometrical invariants of our manifold, for more details of this the reader may consult \cite{MG}.
	
	This cubic form has been important because its profound connections with mirror symmetry and mathematical physics. By encodingc deformations of Calabi-Yau manifolds, the cubic form becomes a central tool for computing mirror maps, predicting enumerative invariants, and determining properties of complex moduli spaces as we can see in \cite{PG}. From a Hodge theoretic perspective the reader can see \cite{DM}.

	In order to extend these powerful algebraic geometric techniques beyond the compact setting, Asakura and Saito investigated in \cite{AS} generalized Jacobian rings for open complete intersections. They obtained analogous results about the primitive cohomology of the logarithmic complex and a generalization of the infinitesimal Torelli theorem.
	
	A cubic form for a cubic threefold and a anticanonical K3 surface has been studied in \cite{AgII}. Subsequent studied was carried in \cite{Wang}. However their construction doesnt resemble the construction of the Yukawa-Griffiths, this is mainly due to the loss of symmetry in the Hodge filtration. We remedy this by working with a half-anticanonical pairs, at the price of working with a different complex, see section \ref{dos}.
	
	By using infinitesimal invariants in mixed Hodge theory, Aguilar studied dualities in \cite{AgI} which allowed him to deduced a series of conjecture analogues to those of Clemens about rational curves. It turned out that the setting for this conjecture were half-anticanonical smooth pairs, this is pairs $(X,Y)$ with $X$ a smooth projective threefold and $Y$ a smooth divisor such that $K_X+2Y=0$. The conjectures were later proved by Aguilar-Zahariuc in \cite{AgII}, for the case when $X$ is Fano.
	
	The objective of this paper is to construct a rigorous algebraic geometric framework for studying an special case of this pairs $(X,Y)$, where $X \subset \mathbb{P}^4$ is a smooth cubic hypersurface and $Y = X \cap H$ is a hyperplane section intersecting $X$ transversally and then construct a non-degenerated cubic form for the pair $(X,Y)$ also, using a suitable bundle $\Sigma=T_X(-\text{log}\ Y)(-Y)$, see subsection \ref{hap}, we establish a precise isomorphism between some cohomology groups associeted to the pair $(X,Y)$ and the graded pieces of a suitable quotient of a polynomial ring.
	
	Although cubic forms can be defined in this general setting of half-anticanonical pairs using the bundle $\Sigma$, this abstract procedure does not guarantee the non-degeneracy of the form nor does it allow for explicit computations. For this reason, we adopt a different perspective to approach the problem.
	
	An important fact in our approach is proving that the associated logarithmic quotient ring $R_{\log}$ exhibits the structure of an Artinian-Gorenstein ring. The presence of the Artinian-Gorenstein duality is indispensable, it provides a canonical trace map down to $\mathbb{C}$ and guarantees the existence and non-degeneracy of our cubic form $C \colon H^1(X,\Sigma)^{\otimes 3} \to \mathbb{C}$, that has important consquences about the infinitesimal invariants of the pair $(X,Y)$. The principal result that we will prove in this paper is the following.
	\begin{theorem}
		Let $X\subset\mathbb{P}^4$ be a cubic smoth hypersurface defined by a homogeneous polynomial $F$ and $Y=\{F=0, z_0=0\}$ an hyperplane section wich intersects $X$ transversally. Denote $\Sigma=T_X({-\text{log}\ Y})(-Y)$, $S=\mathbb{C}[z_0,\dots,z_4]$, $J_{log}=\langle F,\partial_1F,\dots,\partial_4F\rangle$ and $R_{log}=S/J_{log}$, then, we have the following isomorphisms
		\[H^p(X,\Omega_X^{3-p}(\text{log}\ Y)(1-p))\cong H^p(X,\Lambda^p\Sigma)\cong R_{log,2p}\]
		where $R_{log}$ is an Artinian-Gorenstein ring. Moreover, there is a non-degenereted cubic form
		\[C:H^1(X,\Sigma)^{\otimes 3}\to \mathbb{C}.\]
	\end{theorem}
	
	This work is organized as follows. Section 1 introduces preliminary concepts, beginning with Griffiths' classical work on compact Calabi-Yau threefolds, we then develop the necessary framework for pairs $(X,Y)$ and specific algebras that encode cohomological data into algebraic terms, concluding with general constructions for half-anticanonical pairs. Section 2 focuses on the specific case of a cubic threefold-hyperplane pair. In this section, we establish isomorphisms between suitable cohomology groups and the graded pieces of an Artinian-Gorenstein algebra, which ultimately allows us to define a non-degenerate cubic form for our pair in the final subsection.
	\subsection*{AI Usage Disclousure}
	During the preparation of this work, the author used generative AI tools only for the purpose of improving language, grammar, and sentence structure to enhance readability. The
	author takes full responsibility for the content of the paper.
	
	\subsection*{Acknowlegments}
	I want to express my gratitude to my master's thesis advisor, Dr. Rodolfo Aguilar. Beyond his remarkable mathematical insight, his guide has been the absolute foundation of this work, which forms an essential part of my master's thesis. Thanks for introducing me to this field with such genuine passion, for encouraging me to write this text, and for always believing that I could pull it off.
	
	I am immensely grateful for his suggestion of this fascinating problem and for pointing me toward the right methods to attack it. His patience, sharp critique, and valuable suggestions turned every questions into a chance to learn. Above all, thank you for your unconditional support and human kindness through every academic problem I faced along the way. Your constant backing is the real reason this project could write, and your mentorship has left an indelible mark on me as a mathematician. Thank you for everything.
	
	\section{Preliminaries}
	Let us recall the classical result studied by Griffiths in \cite{Griffiths} about a smooth projective Calabi-Yau threefold $X$ defined by a homogeneous polynomial $F$ of degree $d$. In this setting, we have a useful relation between the primitive cohomology groups of our manifold and a suitable graded part in this associated Jacobian ring:
	\[H^{3-p,p}(X)_{\text{prim}}\cong R_{pd+d-5}\]
	where $R=\frac{S}{J_F}=\frac{\mathbb{C}[x_0,x_1,x_2,x_3,x_4]}{\langle \partial_0 F,\partial_1 F,\partial_2 F,\partial_3 F,\partial_4 F\rangle}$ is the graded Jacobian ring.
	\vskip3mm
	Moreover, Macaulay's theorem says that $\text{dim}_\mathbb{C}R_{5d-10}=1$, and the multiplication $R_k\times R_{5d-10-k}\to R_{5d-10}$ defines a perfect pairing, which is identified with the Gauss-Manin connection for $k=d$. In this way, the cubic form $H^1(T_X)^{\otimes 3}\to \mathbb{C}$ defined in terms of the differential of the period map can be calculated using simple polynomial multiplication and the fact that $H^1(T_X\otimes K_X)\cong H^1(\Omega_X^2)$. In this setting, the non-degeneracy of this cubic form has important geometric consequences for this manifold and its period map; in particular, it satisfies the infinitesimal Torelli theorem.
	
	\subsection{Logarithmic bundles for pairs $(X,Y)$}
	
	 In this paper, we are going to work with pairs $(X,Y)$ where $Y$ is a normal crossing divisor and hence with the open variety $U=X\setminus Y$, for this, we recall the relation between the middle cohomology of $U$ and the logarithmic bundles along $U$.
	\begin{definition}
		Let $Y\subseteq X$ be a normal crossing divisor. The de Rham logarithmic bundle $\Omega_X^k(\text{log} \ Y)$ is the subbundle of the bundle of meromorphic forms on $X\setminus Y$, $\Omega_X(*Y)$, consisting of the $k$-forms $\alpha$ such that $\alpha$ and $d\alpha$ admit poles of order at most $1$ along $Y$.
		
		The dual of $\Omega_X^1(\text{log}\ Y)$ is the logarithmic tangent bundle $T_X(-\text{log}\ Y)$, which is the subbundle of $T_X$ consisting of the derivations that preserve the vanishing ideal bundle $\mathcal{O}_X(-Y)$.
	\end{definition}
	With these bundles, we can form the logarithmic de Rham complex
	\[\Omega_X^{\bullet}(\text{log}\ Y):\hspace{5mm}\mathcal{O}_X\rightarrow \Omega_X^{1}(\text{log}\ Y)\rightarrow \dots \rightarrow \Omega_X^{n}(\text{log}\ Y)\rightarrow 0,\]
	following Deligne in \cite{DeligneII}, the cohomology of this complex computes the cohomology with complex coefficients of our open variety
	\[H^k(U,\mathbb{C})\cong\mathbb{H}^k(X,\Omega_X^{\bullet}(\text{log}\ Y)).\]
	Thus, the de Rham-Hodge filtration on the logarithmic complex induces a filtration on $H^\bullet(U,\mathbb{C})$ wich associeted spectral sequence degenerates at $E_1$, hence, we have a descomposition
	\[\text{Gr}_F H^{k}(U,\mathbb{C})\cong\bigoplus_{p+q=k}H^q(X,\Omega_X^{p}(\text{log}\ Y)).\]
	For our case, we are particularly interested in the middle cohomology group:
	\[\text{Gr}_F H^3(U,\mathbb{C})\cong\bigoplus_{p+q=3}H^q(X,\Omega_X^{p}(\text{log}\ Y)).\]
	
	\subsection{Artinian-Gorenstein algebras}
	We will need some definitions and algebraic results about some kind of graded algebras.
	\begin{definition}
		A graded $K$-algebra $R=\bigoplus_{i=0}^{+\infty}R_i$ is called a Artinian-Gorenstein algebra provided the following conditions holds:
		\begin{itemize}
			\item There is $\sigma\geq 0$ such that $\text{dim}_KR_\sigma=1$ and $\text{dim}_KR_i=0$ if $i>\sigma$.
			\item For $0\leq i\leq \sigma$, the multiplication induce a perfect paring
			\[R_i\times R_{\sigma-i}\to R_\sigma.\]
		\end{itemize}
		The perfect paring induce isomorphisms $R_i\cong (R_{\sigma-i})^*$, this phenomenom is known as Gorenstein's duality.
	The number $\sigma$, or sometimes the graded part $R_\sigma$, is called the socle of the algebra.
	\end{definition}
	The following result is due to Macaulay and his proof can be revised in \cite{VosinII}.
	\begin{theorem}[Macaulay]\label{Maculay}
		Let $F_0,F_1,\dots,F_n$ be homogeneous non-constant polynomials in $\mathbb{C}[z_0,\dots,z_n]$ with no non-zero common roots, if $d_i$ denotes the degree of $F_i$, then
		\[\frac{\mathbb{C}[z_0,\dots,z_n]}{\langle F_0,\dots,F_n\rangle}\]
		is an Artinian-Gorenstein algebra with socle $\sigma=\sum_{i=0}^{n}(d_i-1)$.
	\end{theorem}
	
	\subsection{Half-anticanoncal pairs}\label{hap}
	We turn our atention in the case of pairs $(X,Y)$ where $X$ is a projective threfold and $K_X+2Y=0$ in the Picard group. The properties of the twisted logarithmic cohomology groups allows to define a cubic form in this setting, that arise anlagously of the classic Yukawa-Griffiths cubic form. 
	
	For this case, we will consider a particular space of deformations to define the cubic, related to the bundle $\Sigma=T_X(-\text{log}\ Y)(-Y)$. In this setting, the isomorphism $K_X(\text{log}\ Y)(Y)\cong \mathcal{O}_X$ and the trace map $\tau:H^3(X,K_X)\to\mathbb{C}$ gives the following useful result.
	\begin{lemma}
		We have isomorphisms $H^0(X,K_X(\text{log}\ Y)(Y))\cong H^3(X,K_X)\cong \mathbb{C}$ and hence
		\[\text{Hom}(H^0(K_X(\text{log}\ Y)(Y)),H^3(X,K_X))\cong \mathbb{C}.\]
	\end{lemma}
	
	 Taking $\gamma_1,\gamma_2,\gamma_3\in H^1(X,\Sigma)$ 
	we consider the contraction maps induced by this elements and obtain a map
	\[H^0(K_X(\text{log}\ Y)(Y))\xrightarrow{\gamma_1}H^1(X,\Omega_X^2(\text{log}\ Y))\xrightarrow{\gamma_2}H^2(X,\Omega_X^1(\text{log}\ Y)(-Y))\xrightarrow{\gamma_3}H^3(X,K_X).\]
	Using the previous lemma and the fact that the order of the contraction by the elements $\gamma_i$ does not affect the final result, we give the following.
	\begin{proposition}
		The before map yields a cubic form $$H^1(X,\Sigma)^{\otimes 3}\to  \text{Hom}(H^0(K_X(\text{log}\ Y)(Y)),H^3(X,K_X))\cong \mathbb{C}.$$
	\end{proposition}
	 Nevertheless, in this general case, the before procedure does not say the nature of this cubic form, in particular, a priori it can be degenerated.
	
	To further illuminate the geometric nature of this abstractly defined cubic form, we introduce a customized class of Jacobian rings in the spirit of Asakura's foundational work. However, while Asakura's general framework on generalized Jacobian rings provides deep structural insights, it is not directly applycable to our specific objectives, as the rings derived from his standard approach do not naturally align with the exact graded pieces of the logarithmic cohomology required to explicitly construct and evaluate our cubic form.
	
	\section{A cubic threefold-hyperplane as a half-anticanonical pair}\label{dos}

	\subsection{The Twisted log Tangent Bundle}
	
	Now we consider a cubic hypersurface $X =\{F=0\} \subset \mathbb{P}^4$ of degree $d=3$ and a hyperplane section $Y =\{F=0, G=0\}$ with $G$ of degree $e=1$. We assume that $F=0$ and $G=0$ intersect transversally, that is, the gradients $\nabla F$ and $\nabla G$ are linearly independent at each point of $Y$. Although many theoretical aspects can be developed in this setteing, for the purpose of explicit computations and clarity, we will focus on the case where $G = z_0$.
	
	Before working with the involucrated cohomology groups, we will prove some useful facts about the relevant bundles.
	\begin{proposition}\label{iden}
		In the Picard group $\text{Pic}(X)$, the following equalities hold:
		\[K_X\cong\mathcal{O}_X(-2),\ \ \mathcal{O}_X(Y)\cong\mathcal{O}_X(1),\ \ K_X+2Y=0.\]
	\end{proposition}
	\begin{proof}
		By the adjunction formula for a hypersurface of degree $3$ in $\mathbb{P}^4$, we have
		\[K_X\cong \mathcal{O}_{\mathbb{P}^4}(-5+3)\mid_X=\mathcal{O}_X(2).\]
		Defining the hyperplane $H=\{[z_0:z_1:z_2:z_3:z_4]\in \mathbb{P}^4\mid z_0=0\}$, since $Y=X\cap H$, then $\mathcal{O}_X(Y)\cong\mathcal{O}_X(1)$, and hence
		\[K_X+2Y=K_X\otimes\mathcal{O}_X(2Y)\cong\mathcal{O}_X(-2)\otimes\mathcal{O}_X(2)\cong \mathcal{O}_X.\]
	\end{proof}
	
	This means, that our setting is a particular half-anticanonical pair case. Now we want to obtain an Euler-type sequence for our logarithmic tangent bundle; for this, recall the classical Euler sequence of the projective space:
	\[0\rightarrow \mathcal{O}_{\mathbb{P}^4}\xrightarrow{\psi} \mathcal{O}_{\mathbb{P}^4}(1)^{\oplus 5}\xrightarrow{\xi} T_{\mathbb{P}^4}\rightarrow 0,\]
	where $\psi$ sends $1$ to $(z_0,z_1,z_2,z_3,z_4)$ and $\xi(f_0,f_1,f_2,f_3,f_4)=\sum_{i=0}^{4}f_i\partial_i$. Since, locally, $T_{\mathbb{P}^4}(-\text{log}\ H)$ is generated by $\{z_0\partial_0,\partial_1,\partial_2,\partial_3,\partial_4\}$, we can restrict the Euler sequence and obtain
	\[0\rightarrow \mathcal{O}_{\mathbb{P}^4}\xrightarrow{\psi'} \mathcal{O}_{\mathbb{P}^4}\oplus\mathcal{O}_{\mathbb{P}^4}(1)^{\oplus 4}\xrightarrow{\xi'} T_{\mathbb{P}^4}(-\text{log}\ H)\rightarrow 0,\]
	where $\psi'$ sends $1$ to $(1,z_1,z_2,z_3,z_4)$ and $\xi'(1,f_1,f_2,f_3,f_4)=z_0\partial_0+\sum_{i=1}^{4}f_i\partial_i$.
	
	\begin{proposition}
		The preceding sequence is exact and splits; hence,
		\[T_{\mathbb{P}^4}(-\textnormal{log}\ H)\cong \mathcal{O}_{\mathbb{P}^4}(1)^{\oplus 4}.\]
	\end{proposition}
	\begin{proof}
		The injectivity of $\psi'$ and the surjectivity of $\xi'$ follow immediately from the exactness of the original Euler sequence. Now, note that $(c,f_1,f_2,f_3,f_4)\in \text{Ker}\ \xi'$ if and only if $cz_0\partial_0+\sum_{i=1}^{4}f_i\partial_i$ is a multiple of the Euler field, that is, if and only if $(c,f_1,f_2,f_3,f_4)$ is a multiple of $(1,z_1,z_2,z_3,z_4)=\psi'(1)$; hence $\text{Ker}\ \xi'=\text{Im}\ \psi'$.
		
		Moreover, we can define $s:\mathcal{O}_{\mathbb{P}^4}\oplus\mathcal{O}_{\mathbb{P}^4}(1)^{\oplus 4}\to \mathcal{O}_{\mathbb{P}^4}$ by $s(c,f_1,f_2,f_3,f_4)=c$, which is clearly a section of the sequence. Then the sequence splits:
		\[\mathcal{O}_{\mathbb{P}^4}\oplus T_{\mathbb{P}^4}(-\text{log}\ H)\cong \mathcal{O}_{\mathbb{P}^4}\oplus\mathcal{O}_{\mathbb{P}^4}(1)^{\oplus 4},\]
		and hence
		\[T_{\mathbb{P}^4}(-\text{log}\ H)\cong \mathcal{O}_{\mathbb{P}^4}(1)^{\oplus 4}.\]
	\end{proof}

	Now, we are going to use the normal bundle sequence and the previous results to obtain the bundle $\Sigma$ introduced in the before section and presented him as a kernel in an exact sequence.
	\begin{proposition}
		There is an exact sequence:
		\begin{equation}\label{prolon}
			0\rightarrow T_X(-\textnormal{log}\ Y)(-Y)\rightarrow\mathcal{O}_X^{\oplus 4}\rightarrow \mathcal{O}_X(2)\rightarrow 0.
		\end{equation}
	\end{proposition}
	\begin{proof}
		The normal bundle sequence for the hypersurface $X\subset \mathbb{P}^4$ and the bundle $T_{\mathbb{P}^4}(-\text{log}\ H)$ is
		\[0\rightarrow T_X(-\text{log}\ Y)\rightarrow T_{\mathbb{P}^4}(-\text{log}\ H)|_X\rightarrow \mathcal{N}_{X/\mathbb{P}^4}\rightarrow 0.\]
		Since $\mathcal{N}_{X/\mathbb{P}^4}\cong \mathcal{O}_X(3)$, using Proposition 2.2 we get
		\[0\rightarrow T_X(-\text{log}\ Y)\rightarrow \mathcal{O}_X(1)^{\oplus 4}\rightarrow \mathcal{O}_X(3)\rightarrow 0.\]
		Tensoring this sequence by $\mathcal{O}_X(-Y) \cong \mathcal{O}_X(-1)$ and using the previous relations, we obtain
		\[0\rightarrow T_X(-\text{log}\ Y)(-Y)\rightarrow \mathcal{O}_X^{\oplus 4}\rightarrow \mathcal{O}_X(2)\rightarrow 0.\]
	\end{proof}
	We denote $\Sigma = T_X(-\text{log}\ Y)(-Y) \cong T_X(-\text{log}\ Y)(-1)$ and call it the \textit{twisted log tangent bundle}. In local coordinates of $X$, $\zeta_0,\zeta_1,\zeta_2$ in wich $Y$ is defined locally by $\zeta_0=0$, $\Sigma$ is locally generated by $$\left\{\zeta_0^2\frac{\partial}{\partial \zeta_0},\zeta_0\frac{\partial}{\partial \zeta_1},\zeta_0\frac{\partial}{\partial \zeta_2}\right\}.$$

	 We now turn to studying the exterior powers of the sequence (\ref{prolon}). Recall that for every exact sequence of vector bundles $0\to E\to F\to L\to 0$, where $L$ is a line bundle, and for any $p\geq1$, we have the exact sequence:
	\[0\rightarrow \Lambda^p E\rightarrow\Lambda^p F\rightarrow \Lambda^{p-1}E\otimes L\rightarrow 0.\]
	Applying this to the sequence (\ref{prolon}), we get
	\begin{equation}\label{prolong}
		0\rightarrow \Lambda^p T_X(-\text{log}\ Y)(-Y)\rightarrow\Lambda^p \mathcal{O}_X^{\oplus 4}\rightarrow \Lambda^{p-1} T_X(-\text{log}\ Y)(-Y)\otimes \mathcal{O}_X(2)\rightarrow 0.
	\end{equation}
	Let us study the cohomology groups of the induced long exact sequence. For $p=2,3$, we have:
	
	\begin{align*}
		\dots\to H^{p-1}(X,\Lambda^p \Sigma)\to H^{p-1}(X,\Lambda^p\mathcal{O}_X^{\oplus 4})\to H^{p-1}(X,\Lambda^{p-1}\Sigma\otimes\mathcal{O}_X(2))\to H^p(X,&\Lambda^p \Sigma)\\
		\downarrow &\\
		\dots\leftarrow H^p(X,\Lambda^{p-1}\Sigma\otimes\mathcal{O}_X(2)) \leftarrow H^p(X,\Lambda^p&\mathcal{O}_X^{\oplus 4})
	\end{align*}
	
	Since $\Lambda^p\mathcal{O}_X^{\oplus 4} \cong \mathcal{O}_X^{\oplus \binom{4}{p}}$, the Kodaira vanishing theorem gives $H^{p-1}(X,\mathcal{O}_X^{\oplus \binom{4}{p}}) = H^{p}(X,\mathcal{O}_X^{\oplus \binom{4}{p}}) = 0$. Therefore, the exact sequence yields the isomorphism:
	\begin{equation}\label{isom}
		H^p(X,\Lambda^p\Sigma)\cong H^{p-1}(X,\Lambda^{p-1}\Sigma\otimes\mathcal{O}_X(2)).
	\end{equation}
	On the other hand, for $p\in\{1,2,3\}$, tensoring the p sequence (\ref{prolong}) by $\mathcal{O}_X(2p-2)=\mathcal{O}_X((2p-2)Y)$, we get
	\[
	0\rightarrow \Sigma\otimes\mathcal{O}_X(2p-2)\rightarrow\mathcal{O}_X(2p-2)^{\oplus 4}\rightarrow \mathcal{O}_X(2p)\rightarrow 0.
	\]
	Writing $\Sigma\otimes\mathcal{O}_X(l)=\Sigma(l)$, the first terms of its induced long exact sequence are:
	\begin{equation}
		0\to H^0(X, \Sigma(2p-2))\to H^{0}(X,\mathcal{O}_X(2p-2)^{\oplus 4})\xrightarrow{\phi} H^0(X,\mathcal{O}_X(2p))\rightarrow H^1(X,\Sigma(2p-2))\rightarrow 0.
	\end{equation}
	Let $S=\mathbb{C}[z_1,z_2,z_3,z_4]$, then $H^0(X,\mathcal{O}_X(l))\cong \frac{S_{l}}{\langle F\rangle \cap S_{l}}$ consists of the homogeneous polynomials of degree $l$ in four variables modulo $F$, and $\phi$ is induced by the homomorphism $\varphi: S_{2p-2}^{\oplus 4} \to S_{2p}$ given by
	\[\varphi(a_1,a_2,a_3,a_4)=\sum_{i=1}^4a_i\partial_iF.\]
	Thus, we have presented $H^1(X,\Sigma(2p-2))$ as the cokernel of $\phi$, leading to the following theorem:
	\begin{theorem}
		Let $S=\mathbb{C}[z_1,z_2,z_3,z_4]$ and $J_{\log}=\langle F, \partial_1 F,\partial_2 F,\partial_3 F,\partial_4 F\rangle$. For $p=1,2,3$, we have the following isomorphism:
		\begin{equation}
			H^p(X,\Lambda^p\Sigma)\cong H^1(X,\Sigma(2p-2))\cong \operatorname{coker} \phi \cong \left(\frac{S}{J_{\log}}\right)_{2p}.
		\end{equation}
	\end{theorem}
	\begin{proof}
		For $p=1$, the first isomorphism is an identity, while isomorphism (\ref{isom}) provides it for $p=2,3$. The preceding discussion gives the second isomorphism. Moreover,
		\[\operatorname{coker} \phi = \frac{\frac{S_{2p}}{\langle F\rangle\cap S_{2p}}}{\operatorname{Im} \phi},\]
		where 
		\[\operatorname{Im} \phi = \frac{F S_{2p-3} + \langle\partial_1F,\partial_2F,\partial_3F,\partial_4F\rangle_{2p}}{\langle F\rangle \cap S_{2p}},\]
		hence
		\[\operatorname{coker} \phi \cong \frac{S_{2p}}{J_{\log,2p}}=\left(\frac{S}{J_{\log}}\right)_{2p}.\]
		This completes the proof.
	\end{proof}
	\begin{remark}
	Note that if $P\in \mathbb{C}^5$ is a common zero of $F,\partial_1F,\partial_2F,\partial_3F,\partial_4F$, then, the Euler identity tells
	\[z_0\partial_0F(P)=3F(P)-\sum_{i=1}^4z_i\partial_iF(P)=0,\]
	if $z_0(P)=P_0=0$, then $P\in Y$ and $\triangledown F(P)$ is a scalar multiply of $\triangledown z_0(P)$, that contradicts the assumption of $X$ and $H$ interscts transversally, therefore $\partial_0F(P)=0$ and hence $P=0$, since $X$ is smooth.
	\end{remark}
	This yields the following corollary:
	\begin{corollary}
		Let $S$ and $J_{log}$ as in the preceding theorem. Then $\frac{S}{J_{log}}$ is an Artinian-Gorenstein algebra with socle $\sigma=6$.
	\end{corollary}
	\begin{proof}
		The previous remark tells that $F,\partial_1F,\partial_2F,\partial_3F,\partial_4F$ satisfy the conditions of Macaulay's theorem (\ref{Maculay}). This gives the desired.
	\end{proof}
	
	\subsection{The logarithmic Jacobian ring}
	Now, is turn to relate the previous results with the cohomology of our open variety. We have a canonical isomorphism
	\[\Lambda^p\Sigma\cong(\Lambda^{3-p}\Sigma)^*\otimes\text{det}(\Sigma),\]
	where $\text{det}(\Sigma)=\Lambda^3 T_X(-\text{log} Y)(-Y)\cong \text{det}(\Omega_X^1(\text{log}\ Y)\otimes \mathcal{O}_X(Y))^*\cong \text{det}(\Omega_X^1(\text{log}\ Y))^*\otimes \mathcal{O}_X(-Y)^{\otimes 3}$,
	and $\text{det}(\Omega_X^1(\text{log}\ Y))=K_X\otimes\mathcal{O}_X(Y)$, using the identities of (\ref{iden}) we obtain
	\[\text{det}(\Sigma)\cong (K_X\otimes \mathcal{O}_X(1))^*\otimes\mathcal{O}_X(-3)\cong K_X\otimes\mathcal{O}_X(-4)\cong \mathcal{O}_X(-2)\cong K_X,\]
	this yields the following bundle isomorphism
	 \[\Lambda^p\Sigma\cong\Lambda^{3-p}\Omega_X^1(\text{log}\  Y)(Y)\otimes K_X\cong\Omega_X^{3-p}(\text{log}\ Y)(3-p)\otimes\mathcal{O}_X(2)=\Omega_X^{3-p}(\text{log}\ Y)(1-p).\]
	 All the preceding results allow us to stablish the following theorem.
	 \begin{theorem}\label{Goren}
	 	Let $F\in S=\mathbb{C}[z_0,z_1,z_2,z_3,z_4]$ be an homogeneous polynomial of degree $d=3$ and $X=\{F=0\}\subset \mathbb{P}^4$ the hypersurface defined by $F$. If $X$ is smooth and the hiperplane $H=\{z_0=0\}\subset \mathbb{P}^4$ intersects $X$ transversally, then, for every $p=0,1,2,3$ the following isomorphisms holds
	 	\begin{equation}
	 		H^{p}(X,\Omega_X^{3-p}(\textnormal{log}\ Y)(1-p))\cong H^p(X,\Lambda^p\Sigma)\cong \left(\frac{S}{J_{log}}\right)_{2p},
	 	\end{equation}
	 	where $J_{\log}=\langle F, \partial_1 F,\partial_2 F,\partial_3 F,\partial_4 F\rangle$ and there is a perfect paring
	 	\begin{equation}
	 		H^{p}(X,\Omega_X^{3-p}(\textnormal{log}\ Y)(1-p))\times H^{3-p}(X,\Omega_X^{p}(\textnormal{log}\ Y)(p-2))\to H^{3}(X,\Omega_X^{0}(\textnormal{log}\ Y)(-2))\cong \mathbb{C}.
	 	\end{equation}
	 	\begin{proof}
	 		For the isomorphisms, only remains the case $p=0$ that follows from
	 		\[\Omega_X^3(\text{log}\ Y)(1)\cong K_X(2)\cong \mathcal{O}_X\]
	 		and $H^0(X.\mathcal{O}_X)\cong\mathbb{C}=\left(S/J_{log}\right)_0$.
	 		The perfect paring is induced by the perfect paring in the Artinian-Gorenstein algebra $\frac{S}{J_{log}}$, with the consideration $ H^3(X,\Omega_X^0(\text{log}\ Y)(-2))\cong \left(\frac{S}{J_{log}}\right)_6\cong\mathbb{C}$.
	 	\end{proof}
	 \end{theorem}
	
	We will denote $S/J_{log}$ as $R_{log}$ and call it the logarithmic Jacobian ring.
	\begin{corollary}
		We have calculated the middle cohomology group of $U = X \setminus Y$ as
		\[ H^1(X,\Omega_X^2(\text{log}\ Y)) \cong R_{\log, 2}, \quad H^2(X, \Omega_X^1(\log Y)) \cong R_{\log, 5}\] \[H^0(X, \Omega_X^3(\log Y)) \cong  H^3(X, \Omega_X^0(\log Y)) \cong 0. \]
	\end{corollary}
	\begin{proof}
		The first isomorphism follows from the preceding theorem by setting $p=1$. For the second, Serre duality gives
		\[ H^2(X, \Omega_X^1(\log Y)) \cong H^1(X, T_X(-\log Y) \otimes K_X)^*. \]
		Since $T_X(-\log Y) \otimes K_X \cong \Sigma(-Y)$, tensoring the short exact sequence (\ref{prolon}) by $\mathcal{O}_X(-Y)$ and passing to the long exact sequence in cohomology, we obtain
		\[ H^1(X, \Sigma(-Y)) \cong H^0(X, \mathcal{O}_X(-1)) \cong S_1. \]
		Because the generators of $J_{\log}$ have degree at least $2$, we have $S_1 \cong R_{\log, 1}$, and thus
		\[ H^2(X, \Omega_X^1(\log Y)) \cong H^1(X, \Sigma(-Y))^* \cong (R_{\log, 1})^* \cong R_{\log, 5}. \]
		For the last two isomorphisms, we note that
		\[ H^3(X, \Omega_X^0(\log Y)) \cong H^3(X, \mathcal{O}_X) = 0, \]
		\[ H^0(X, \Omega_X^3(\log Y)) \cong H^0(X, K_X(Y)) \cong H^0(X, \mathcal{O}_X(-1)) = 0. \]
	\end{proof}
	
	\begin{proposition}
		The homogeneous components of the logarithmic Jacobian ring $R_{\log}$ have the following dimensions:
		\[ \dim_\mathbb{C}(R_{\log, 0}) = \dim_\mathbb{C}(R_{\log, 6}) = 1, \]
		\[ \dim_\mathbb{C}(R_{\log, 1}) = \dim_\mathbb{C}(R_{\log, 5}) = 5, \]
		\[ \dim_\mathbb{C}(R_{\log, 2}) = \dim_\mathbb{C}(R_{\log, 4}) = 11. \]
	\end{proposition}
	\begin{proof}
		Gorenstein duality implies that it suffices to compute the dimensions for degrees $0$, $1$, and $2$. The first is trivial since $S_0 = \mathbb{C}$. The second follows from the fact that $J_{\log}$ contains only polynomials of degree at least $2$. The third follows from the quotient structure
		\[ R_{\log, 2} \cong \frac{S_2}{\langle \partial_1 F, \partial_2 F, \partial_3 F, \partial_4 F \rangle_2}, \]
		which yields
		\[ \dim_\mathbb{C}(R_{\log, 2}) = \dim_\mathbb{C}(S_2) - \dim_\mathbb{C}(\langle \partial_1 F, \partial_2 F, \partial_3 F, \partial_4 F \rangle_2) = \binom{6}{4} - 4 = 11. \]
	\end{proof}
	
	\begin{corollary}
		The non-trivial Hodge numbers of the open variety $U$ are
		\[ h^{2,1} = 11, \quad h^{1,2} = 5, \quad h^{3,0} = h^{0,3} = 0.\]
	\end{corollary}
	
	\subsection{A Cubic Form}
	In \cite{AGG}, R. Aguilar, M. Green and P. Griffiths constructed cubic forms for a pair $(X,Y)$ with $X$ a cubic threefold in $\mathbb{P}^4$ and $Y$ an anticanonical smoth section defined by a quadric polinomial and use this to prove a version of the infinitesimal Torelli theorem. Also, Z. Wang prove some other properties of this cubic form and used another technicals for give infinitesimal invariants in this kind of pairs in \cite{Wang}.
	
	Inspired in these works, for our special case, we construct a cubic form on $H^1(X, \Sigma)$ a little more similar to the classic case, using the perfect pairing of Theorem \ref{Goren} and a suitable trace map over $H^3(X, K_X)$.
	
	\begin{proposition}
		The space $H^3(X, K_X)$ is isomorphic to $R_{\log, 6}$, and there exists an isomorphism called the trace map:
		\[ \tau : H^3(X, K_X) \to \mathbb{C}. \]
	\end{proposition}
	\begin{proof}
		Theorem~\ref{Goren} yields an isomorphism $\varphi' : R_{\log, 6} \to H^3(X, \Omega_X^0(\log Y)(-2))$. Since $\Omega_X^0(\log Y)(-2) \cong \mathcal{O}_X(-2) \cong K_X$, this induces an isomorphism
		\[ \varphi : R_{\log, 6} \to H^3(X, K_X). \]
		Using integration, we obtain the trace map
		\[ \tau : H^3(X, K_X) \to \mathbb{C}. \]
		Since this map is non-zero and both spaces are $1$-dimensional, $\tau$ is an isomorphism.
	\end{proof}
	
	For integers $0 \leq i, j \leq 6$ such that $i+j \leq 6$ and for any $[P] \in R_{\log, i}$, polynomial multiplication defines a map $R_{\log, j} \to R_{\log, i+j}$. This induces a morphism
	\[ \mu : R_{\log, i} \to \operatorname{Hom}(R_{\log, j}, R_{\log, i+j}). \]
	
	\begin{proposition}\label{inj}
		For integers $0 \leq i, j \leq 6$ such that $i+j \leq 6$, the preceding morphism is injective.
	\end{proposition}
	\begin{proof}
		If $[P] \in \ker \mu$, then $[P \cdot A] = 0$ for every $[A] \in R_{\log, j}$. Consequently, $[P \cdot A \cdot B] = 0$ for all $[B] \in R_{\log, 6-i-j}$. However, the multiplication map
		\[ R_{\log, 6-i-j} \times R_{\log, j} \to R_{\log, 6-i}, \quad ([A], [B]) \mapsto [A \cdot B] \]
		is surjective; thus, $[P \cdot Q] = 0$ for every $[Q]\in R_{\log, 6-i}$. Macaulay's theorem then implies that $[P]=0$.
	\end{proof}
	
	Consider the trilinear multiplication map
	\[ M : R_{\log, 2} \times R_{\log, 2} \times R_{\log, 2} \to R_{\log, 6}. \]
	Suppose that $[P_1] \in R_{\log, 2}$ satisfies $[P_1 P_2 P_3] = 0$ for every $[P_2], [P_3] \in R_{\log, 2}$. Then $[P_1] \in \ker \mu$, and hence $[P_1] = 0$. This shows that the map is non-degenerate in the first factor. Similar arguments demonstrate that it is non-degenerate in each factor.
	
	Moreover, Theorem~\ref{isom} gives an isomorphism $\psi : H^1(X, \Sigma) \xrightarrow{\sim} R_{\log, 2}$, which leads to the following result.
	
	\begin{theorem}
		The map $C' = \tau \circ \varphi \circ M \circ \psi^{\times 3} : H^1(X, \Sigma)^{\times 3} \to \mathbb{C}$ induces a non-degenerate cubic form
		\[ C : H^1(X, \Sigma)^{\otimes 3} \to \mathbb{C}. \]
	\end{theorem}
	\begin{proof}
		The multilinearity of $M$ allows us to define the cubic form $C$. Since $\tau$, $\varphi$, and $\psi$ are isomorphisms and $M$ is non-degenerate, it follows that $C$ is non-degenerate.
	\end{proof}

	\bibliographystyle{abbrv}
	\bibliography{bib.bib}
	\nocite{*}
	
	\vspace{2em} \noindent
	\textsc{Gustavo Serrano Herrera} \\
	\small\textsc{Centro de Investigación en Matemáticas, A.C. (CIMAT)} \\
	\small\textsc{Jalisco S/N, Col. Valenciana, C.P. 36023, Guanajuato, Gto., México} \\
	\small\text{E-mail address}: \href{mailto:gserranoh001@gmail.com}{\texttt{gserranoh001@gmail.com}}
\end{document}